\documentclass[12pt]{article}   	
\usepackage[margin=1.3in]{geometry}                		
\usepackage{graphicx}				
\usepackage{amssymb}

\usepackage{hyperref}\hypersetup{colorlinks,
  linkcolor={blue},
  citecolor={blue},
  urlcolor=blue
}

\def\arXiv#1{{\href{http://front.math.ucdavis.edu/#1}{arXiv:\linebreak[0]#1}}}

\begin{document}

\title{Equivariant Characteristic Classes\\
(Commentary on [116])}
\author{Loring W.\ Tu}
\date{}							

\maketitle
\begin{abstract}
This is a commentary on Raoul Bott and Loring Tu's joint article "Equivariant characteristic classes in the Cartan model," which appeared in "Geometry, Analysis, and Applications (Varanasi, 2000)," World Scientif Publishing, River Edge, NJ, 3--20. The article is also included in "Raoul Bott: Collected Papers," Vol. 5, as article [116]. The commentary discusses the genesis of the article, its influence, and the current state of the problem concerning equivariant characteristic classes.
\end{abstract}

I was trained as an algebraic geometer under Phillip A.\ Griffiths, but I have always had an abiding interest in topology, especially Raoul Bott's kind of topology.  
In 1995 Raoul Bott gave a series of lectures at Brown University on equivariant cohomology.
I was very much captivated by his presentation of the subject matter.
I decided to reorient myself and work with Raoul Bott on equivariant
cohomology.
Our collaboration had the advantage of geographical proximity---I was
teaching at Tufts University, only two subway stops from Harvard,
where Raoul was.

In one of his courses he showed the students a way of computing the equivariant cohomology of the projective space under a circle action.
He suggested to me the problem of computing the equivariant cohomology
of the complete flag manifold $U(n)/T$, where $T =  U(1) \times \cdots
\times U(1)$
is the maximal torus in the unitary group $U(n)$, under
the left action of $T$.  He
had a conjectural formula not only for the this case, but more
generally for the homogeneous space $G/H$,
where $G$ is a compact connected Lie group and $H$ is a closed
subgroup of maximal rank.
(The \emph{rank} of a compact Lie group is the dimension of the maximal torus it contains.)
Using the method Bott showed us in class, I was able to prove his
conjecture. 

In my excitement, I suggested to him that we could publish this as a joint paper.
Alas, it was not to be, for it turned out that the equivariant
cohomology of such a homogeneous space under the maximal torus action had already been worked out by
Alberto Arabia (\cite{arabia86}, \cite{arabia89}, \cite{brion}).
Although Bott's method was different and original, it was not enough for a paper.

Perhaps to assuage my disappointment, Raoul then suggested to me a problem that he was sure no one had worked out yet. 
 
There are two approaches to defining the characteristic classes of a 
vector bundle or a principal $K$-bundle for a Lie group $K$ in the smooth category.
The first is topological, as elements of the cohomology ring of the classifying space $BK$ of the principal $K$-bundle.
The second is differential-geometric, as certain differential forms constructed from the curvature of a connection on the bundle.
Since the cohomology ring $H^*(BK)$ consists of invariant polynomials on the Lie algebra of $K$, the Chern--Weil homomorphism in fact relates the topological approach to the differential-geometric approach.

For a $G$-equivariant principal $K$-bundle $\pi\colon P \to M$,
again there are two approaches to defining equivariant characteristic classes.
In the topological approach, one forms the homotopy quotients $P_G$ and $M_G$ and defines the equivariant characteristic classes of $P \to M$ as the ordinary characteristic classes of the principal $K$-bundle $P_G \to M_G$.
In the differential-geometric approach, equivariant characteristic classes are equivariant differential forms constructed from the curvature of a connection on the bundle and a map associated to the action of $G$ sometimes called the moment map.

Raoul told me that the topological approach was popular in the United
States, but the differential-geometric approach dominated in France,
as in the works of Nicole Berline  and Mich\`ele Vergne
(\cite{berline--vergne82}, \cite{berline--vergne83}).
It was generally assumed that the two approaches are equivalent---they describe the same objects and yield the same theorems, but no one had actually shown explicitly the equivalence.

Raoul had by then retired from Harvard and was spending most of the winter months at the University of California in San Diego, near his daughter Jocelyn's home in Rancho Santa Fe.
Since my parents and my brother Charles lived near San Diego
and Charles was the Associate Dean of the Engineering School at UCSD,
it was quite convenient for me to meet Raoul at UCSD during my winter
breaks, three thousand miles away from my home institution of Tufts.
I still remember discussing the problem of equivariant characteristic classes with him in his office at UCSD.
He outlined to me a way to bridge the two approaches for a circle action.
I got very excited, because it was such a beautiful and unexpected construction.
He asked me to generalize it to a torus action.
Our discussion ended there, for he wanted to go home before the terrible rush hour traffic in Southern California.
He asked me if I would like to see his daughter's house.
So I followed him in my car on a narrow, winding road along the Pacific Coast from La Jolla to Rancho Santa Fe.
It was a magnificent family estate situated on a hill with a panoramic view of the ocean.  Raoul and his wife Phyllis stayed in the large guest house on the property.
With such an abode in the temperate climate of Southern California,
no wonder they eventually moved there.

After returning to Cambridge, I worked out the generalization of Raoul's construction to a torus action.  
By then he had returned to Harvard.
When I showed it to him, he said, ``This is very good, but it is not
publishable.  You need to do it for a compact Lie group action.''
This was typical of Raoul, to suggest a problem in a sequence of steps.
I was stuck for a couple of years, I think, before one day I
suddenly saw how to generalize it to a connected compact Lie group
action.
We agreed to meet after he attended church service at Saint Paul's in
Cambridge one Sunday, because the church was close to my apartment.
We met in an Au Bon Pain near the church and I outlined my
solution to him.
We decided to make it a joint paper.

At the time, for me, this was a technical exercise.
Only much later did I realize that we had discovered the equivariant
analogue of the Chern-Weil homomorphism.
At the time Raoul was asked by
some Indian conference to contribute a paper, so we gave it to them.
It was probably one of my best papers, but it ended up in an obscure
Indian conference proceeding that few people have access to.
Fortunately, the paper has been on ArXiv and has found an audience
there.

In 2004 the Botts moved permanently to California.
I continued to meet with Raoul whenever I visited my parents.
Each time he would suggest some interesting problem for me to work on,
but [116] is our last successful research collaboration.
For one reason or another, I was never able to carry the other
projects to fruition.
His generosity of spirit, patience, and encouragement and his
inimitable
lecture style remain for me a model to emulate.
His passing in December 2005 was for me a profound loss and a source
of great sadness.

I am aware of two works that are in some sense inspired by [116].
A Lie groupoid $G$ is an object with all the properties of a Lie group
except that the multiplication map is defined only on a subset of
$G \times G$.
In \cite{laurent-gengoux--tu--xu}  Laurent-Gengoux, Tu, and Xu
generalized our work to principal 
bundles over a Lie groupoid.
In [116] Bott and I assumed that the Lie group acting on the spaces is
compact and connected.
The connectedness hypothesis is in fact not necessary.
Using a completely different approach, Andreas K\"ubel and Andreas
Thom \cite{kubel--thom} reprove our theorem and
remove the connectedness hypothesis.
To deal with the infinite-dimensionality of the universal bundle of a
Lie group, Bott and I approximated the infinite-dimensional spaces
with finite-dimensional manifolds.
K\"ubel and Thom instead represent the infinite-dimensional spaces
as the geometric realization of a semi-simplicial manifold
and use semi-simplicial de Rham theory in place of de Rham theory.
Since the Cartan model is valid only for a compact Lie group action,
the work of K\"ubel and Thom completes the story for such an
                 action.
The equivariant cohomology of a noncompact Lie group action, however,
remains a
largely unexplored virgin territory.

\bibliographystyle{alpha}

\bigskip
\noindent
{\sc Department of Mathematics, Tufts University, Medford, MA 02155}

\noindent
\textit{Email address}:  \texttt{loring.tu@tufts.edu}

\end{document}